


\font\eightpt=cmr8
\font\ninept=cmr9


\font\mathbold=cmmib10

\font\ninerm=cmr9     \font\eightrm=cmr8   \font\sixrm=cmr6      
\font\ninei=cmmi9     \font\eighti=cmmi8   \font\sixi=cmmi6      
\font\ninesy=cmsy9    \font\eightsy=cmsy8  \font\sixsy=cmsy6     
\font\ninebf=cmbx9    \font\eightbf=cmbx8  \font\sixbf=cmbx6     
\font\ninett=cmtt9    \font\eighttt=cmtt8                        
\font\nineit=cmti9    \font\eightit=cmti8     
\font\ninesl=cmsl9    \font\eightsl=cmsl8                        

\font\tensc=cmcsc10   \font\ninesc=cmcsc9  \font\eightsc=cmcsc8  

\font\eightssq=cmssq8  \font\eightssqi=cmssqi8  


\font\tenssbx=cmssbx10 

  \font\twelvebf=cmbx12
  
\def\sc{\tensc}  \def\mc{\ninerm}


\newskip\ttglue
\def\tenpoint{\def\rm{\fam0\tenrm}%
  \textfont0=\tenrm \scriptfont0=\sevenrm \scriptscriptfont0=\fiverm
  \textfont1=\teni  \scriptfont1=\seveni  \scriptscriptfont1=\fivei
  \textfont2=\tensy \scriptfont2=\sevensy \scriptscriptfont2=\fivesy
  \textfont3=\tenex \scriptfont3=\tenex   \scriptscriptfont3=\tenex
  \textfont\itfam=\tenit  \def\it{\fam\itfam\tenit}%
  \textfont\slfam=\tensl  \def\sl{\fam\slfam\tensl}%
  \textfont\ttfam=\tentt  \def\tt{\fam\ttfam\tentt}%
  \textfont\bffam=\tenbf  \scriptfont\bffam=\sevenbf
   \scriptscriptfont\bffam=\fivebf \def\bf{\fam\bffam\tenbf}%
  \tt \ttglue=.5em plus.25em minus.15em
  \normalbaselineskip=12pt
  \setbox\strutbox=\hbox{\vrule height8.5pt depth3.5pt width0pt}%
  \let\sc=\tensc \let\mc=\ninerm  
  \def\cyr{\tencyr\cyracc}\def\cyri{\tencyri\cyracc}
  \let\big=\tenbig  \normalbaselines\rm}

\def\ninepoint{\def\rm{\fam0\ninerm}%
\textfont0=\ninerm  \scriptfont0=\sixrm  \scriptscriptfont0=\fiverm
\textfont1=\ninei   \scriptfont1=\sixi   \scriptscriptfont1=\fivei
\textfont2=\ninesy  \scriptfont2=\sixsy  \scriptscriptfont2=\fivesy
\textfont3=\tenex   \scriptfont3=\tenex  \scriptscriptfont3=\tenex
\textfont\itfam=\nineit  \def\it{\fam\itfam\nineit}%
\textfont\slfam=\ninesl  \def\sl{\fam\slfam\ninesl}%
\textfont\ttfam=\ninett  \def\tt{\fam\ttfam\ninett}%
\textfont\bffam=\ninebf  \scriptfont\bffam=\sixbf
\scriptscriptfont\bffam=\fivebf\def\bf{\fam\bffam\ninebf}%
\tt\ttglue=.5em plus.25em minus.15em
\normalbaselineskip=11pt
\setbox\strutbox=\hbox{\vrule height8pt depth3pt width0pt}%
\let\sc=\ninesc\let\mc=\eightrm
\def\cyr{\ninecyr\cyracc}\def\cyri{\ninecyri\cyracc}
\let\big=\ninebig\normalbaselines\rm}

\def\eightpoint{\def\rm{\fam0\eightrm}%
  \textfont0=\eightrm \scriptfont0=\sixrm \scriptscriptfont0=\fiverm
  \textfont1=\eighti  \scriptfont1=\sixi  \scriptscriptfont1=\fivei
  \textfont2=\eightsy \scriptfont2=\sixsy \scriptscriptfont2=\fivesy
  \textfont3=\tenex   \scriptfont3=\tenex \scriptscriptfont3=\tenex
  \textfont\itfam=\eightit  \def\it{\fam\itfam\eightit}%
  \textfont\slfam=\eightsl  \def\sl{\fam\slfam\eightsl}%
  \textfont\ttfam=\eighttt  \def\tt{\fam\ttfam\eighttt}%
  \textfont\bffam=\eightbf  \scriptfont\bffam=\sixbf
  \normalbaselineskip=9pt
  \let\sc=\eightsc \let\mc=\sevenrm  
  \def\cyr{\eightcyr\cyracc}\def\cyri{\eightcyri\cyracc}
  \let\big=\eightbig  \normalbaselines\rm}%
\def\nospace{\nulldelimiterspace0pt\mathsurround0pt}%
\def\tenbig#1{{\hbox{$\left#1\vbox to8.5pt{}\right.\nospace$}}}%
\def\ninebig#1{{\hbox{$\textfont0=\tenrm\textfont2=\tensy
  \left#1\vbox to7.25pt{}\right.\nospace$}}}%
\def\eightbig#1{{\hbox{$\textfont0=\ninerm\textfont2=\ninesy
  \left#1\vbox to6.5pt{}\right.\nospace$}}}%

\def\nonextendedbold{
  \font\fiveb=cmb10 at 5pt
  \font\sixb=cmb10 at 6pt
  \font\sevenb=cmb10 at 7pt
  \font\eightb=cmb10 at 8pt
  \font\nineb=cmb10 at 9pt
  \font\tenb=cmb10
  \font\twelveb=cmb10 at 12pt
  \let\fivebf=\fiveb
  \let\sixbf=\sixb
  \let\sevenbf=\sevenb
  \let\eightbf=\eightb
  \let\ninebf=\nineb
  \let\tenbf=\tenb
  \let\twelvebf=\twelveb
}

\def\leftrighttop#1#2{
  \headline{\ifnum\pageno=1\hfil\else{\ninept #1 \hfil #2}\fi}
}

\def\firstnopagenum{
  \footline{\ifnum\pageno=1 \hfil \else \hfil{\rm \number\pageno}\hfil\fi}
}

\def\maketitle#1#2#3#4{
  \centerline {\titlefont #1}
  \medskip
  \centerline {\eightpt #2}
  \medskip
  \centerline {\tensc #3}
  \medskip
  \centerline {\tensc #4}
  \bigskip
}


\outer\def\floattext#1 #2. #3\par{
  $$
  \vbox{
    \hsize #1 true in
    \noindent{\bf #2.}\enskip #3
  }
  $$
}


\def\lsection#1\par{
  \bigskip\vskip\parskip
  \leftline{\sectionfont#1}\nobreak\medskip\noindent
}

\def\csection#1\par{
  \bigskip\vskip\parskip
  \centerline{\sectionfont#1}\nobreak\medskip\noindent
}

\def\rsection#1\par{
  \bigskip\vskip\parskip
  \rightline{\sectionfont#1}\nobreak\medskip\noindent
}
\def\section{\lsection}

\def\boldlabel#1. {\noindent{\bf #1.\enspace}}
\def\subsection#1. {\medskip\noindent{\bf #1.\enspace}}



\font\tenfrak=eufm10
\font\sevenfrak=eufm7
\font\fivefrak=eufm5
\newfam\frakfam
\textfont\frakfam=\tenfrak
\scriptfont\frakfam=\sevenfrak
\scriptscriptfont\frakfam=\fivefrak

\def\janksc#1#2 {#1{\eightpt#2}}
\def\jankscsp#1#2 {#1{\eightpt#2}\ }
\def\scproclaim#1.#2\par{\noindent\jankscsp #1.\enspace{\it#2\par}}


\def\ref#1{[#1]}

\def\quote{
  \begingroup
    \baselineskip 10pt
    \parfillskip 0pt
    \interlinepenalty 10000
    \leftskip 0pt plus 40pc minus \parindent
    \let\rm=\quoterm\let\sl=\quotesl\everypar{\sl}
    \obeylines
}
\def\author#1(#2){\nobreak\smallskip\rm--- \rm#1\unskip\enspace(#2)\par\endgroup}

\def\titlefont{\twelvebf}
\def\sectionfont{\tenssbx}
\def\quoterm{\eightssq}
\def\quotesl{\eightssqi}


\def\widemargins{
  \magnification=\magstephalf\hoffset=40pt \voffset=28pt
  \hsize=29pc  \vsize=45pc  \maxdepth=2.2pt  \parindent=19pt
}
\def\bookheader#1#2{
  \nopagenumbers
  \def\leftheadline{{\rm\folio}\hfil{\eightpoint#1}\hfil}
  \def\rightheadline{\hfil{\eightpoint#2}\hfil{\rm\folio}}
  \headline{\ifodd\pageno{\ifnum\pageno<2\hfil\else\rightheadline\fi}\else\leftheadline\fi}
}

\tenpoint



\def\xskip{\hskip 7pt plus 3pt minus 4pt}

\def\proof{\medbreak\noindent{\it Proof.}\xskip\ignorespaces}

\def\slug{\quad\hbox{\kern1.5pt\vrule width2.5pt height6pt depth1.5pt\kern1.5pt}\medskip}
\def\noskipslug{\quad\hbox{\kern1.5pt\vrule width2.5pt height6pt depth1.5pt\kern1.5pt}}

\newdimen\algindent
\newif\ifitempar \itempartrue 
\def\algindentset#1{\setbox0\hbox{{\bf #1.\kern.25em}}\algindent=\wd0\relax}
\def\algbegin #1 #2{\algindentset{#21}\alg #1 #2} 
\def\aalgbegin #1 #2{\algindentset{#211}\alg #1 #2} 
\def\alg#1(#2). {\medbreak 
  \noindent{\bf#1}({\it#2\/}).\xskip\ignorespaces}
\def\algstep#1.{\ifitempar\smallskip\noindent\else\itempartrue
  \hskip-\parindent\fi
  \hbox to\algindent{\bf\hfil #1.\kern.25em}%
  \hangindent=\algindent\hangafter=1\ignorespaces}




\def\op#1{\mathop{\hbox{#1}}\nolimits}





\newcount\thmcount  
\thmcount=1
\newcount\sectcount  
\sectcount=1
\newcount\figcount  
\figcount=1
\newcount\eqcount  
\eqcount=1

\def\oldno#1{\eqno({\oldstyle#1})}
\def\refeq#1{({\oldstyle#1})}
\def\adveq{\oldno{\the\eqcount}\global\advance\eqcount by 1}  
\def\advthm{\the\thmcount\global\advance \thmcount by 1}

\def\advsect{\section\the\sectcount\global\advance\sectcount by 1. }

\def\caption#1{\centerline{\ninepoint{\bf Fig.~\the\figcount\global\advance\figcount by 1.\enspace}#1}}

\outer\def\parenproclaim #1 (#2).#3\par{\medbreak
  \noindent{\bf #1}\enspace\rm({\it #2\/}).\nobreak\ignorespaces{\sl #3\par}
  \ifdim\lastskip<\medskipamount \removelastskip\penalty55\medskip\fi}


\newdimen\axiomindent
\def\axiomindentset#1{\setbox0\hbox{{\bf #1.\kern.25em}}\axiomindent=\wd0\relax}
\def\axiom#1. [#2.]{\ifitempar\par\noindent\else\itempartrue
  \hskip-\parindent\fi%
  \hbox to\axiomindent{\bf\hfil #1.\kern.25em}%
  \hangindent=\axiomindent\hangafter=1[{\it #2.}]}

\input eplain
\input epsf


\def\UH{\op{\mc UH}}
\def\UHsmall{\op{\sevenrm UH}}

\def\tr{\op{\rm tr}}

\def\mfield#1#2{\hbox{\tt#1}(#2)} 

\newcount\refcount
\refcount=1
\def\advref{\global\advance \refcount by 1}
\edef\refblieberger{\the\refcount}\advref
\edef\refbodinione{\the\refcount}\advref
\edef\refbodinitwo{\the\refcount}\advref
\edef\refanalyticcombinatorics{\the\refcount}\advref
\edef\refjanson{\the\refcount}\advref
\edef\refjansonkuba{\the\refcount}\advref
\edef\refkemp{\the\refcount}\advref
\edef\refkirschenhofer{\the\refcount}\advref
\edef\refkubaone{\the\refcount}\advref
\edef\refkubatwo{\the\refcount}\advref
\edef\refmorris{\the\refcount}\advref
\edef\refpemantle{\the\refcount}\advref
\edef\refstevanovic{\the\refcount}\advref
\def\ref#1{[#1]}

\widemargins
\bookheader{CHEUNG, DEVROYE, AND GOH}{A NOTE ON PLANE TREES WITH DECREASING LABELS}

\maketitle{A note on plane trees with decreasing labels}{}{Tsun-Ming Cheung, Luc Devroye, {\rm and} Marcel Goh}{\sl School of Computer Science, McGill University}

\floattext 4.5 \ninebf Abstract.
\ninepoint This note derives asymptotic upper and lower bounds for the number of planted plane trees on $n$ nodes assigned labels from the set $\{1,2,\ldots, k\}$ with the restriction that on any path from the root to a leaf, the labels must strictly decrease. We illustrate an application to calculating the largest eigenvalue of the adjacency matrix of a tree.
\smallskip
\noindent\boldlabel Keywords. Planted plane trees, decreasing labels, eigenvalues of trees.
\smallskip
\noindent\boldlabel Mathematics Subject Classification. 05C05, 05C30, 05C50.

\bigskip\bigskip\noindent
{\sc In this note,} we count the number of {\it planted plane trees}
(also sometimes called {\it ordered trees}) on $n$ nodes, each node is given a label from the set $\{1,2,\ldots,k\}$, and labels must strictly decrease on any path from the root to a leaf. Let $G_{n,k}$ denote this count for given positive integers $n$ and $k$.

It is clear that $G_{1,1} = 1$ and $G_{n,1} = 0$ for all $n\ge 2$. Then, letting
$C(n)$ denote the set of all compositions of $n$, we have the recurrence formula
\edef\eqGkrecurrence{\the\thmcount}
$$G_{n,k} = G_{n,k-1} + \sum_{S\in C(n-1)} \prod_{s\in S} G_{s,k-1}\adveq$$
for $k\ge 2$.
The first term, $k$, corresponds to the number of labellings that do not use the label $k$. For the
second term, the root label must be $k$, the children of the root have subtrees with a total of $n-1$ nodes, and each subtree must only use labels from $\{1,2,\ldots,k-1\}$.

The sequence $G_{n,k}$ has two parameters, which suggests the use of a bivariate generating function in its analysis (see, e.g.,~\ref{\refpemantle} for a detailed account of methods related to multivariate generating functions). However we found that the nature of our specific recurrence made it easiest to work with a sequence of single-variable generating functions instead.

For $k\ge 1$, let $G_k(z)$ be the generating function
$$ G_k(z) = \sum_{n=0}^\infty G_{n,k}z^n.$$
Immediately we see that $G_1(z) = z$,
and from the recurrence~\refeq{\eqGkrecurrence}, we have
\edef\eqGkformula{\the\eqcount}
$$G_k(z) = G_{k-1}(z) + {z\over 1-G_{k-1}(z)}\adveq$$
for $k\ge 2$.
For any fixed $k$, the coefficients $G_{n,k}$ of $G_k(z)$ are nonnegative integers, so by Pringsheim's theorem~\ref{\refanalyticcombinatorics}, if the radius of convergence of $G_k(z)$ is $R_k$, then $R_k$ is a singularity of the function $G_k(z)$.
From the formula~\refeq{\eqGkformula},
it is clear that any pole of $G_{k-1}(z)$ as well as any
solution $z$ to $G_{k-1}(z) = 1$ is a pole of $G_k(z)$.

For convenience, we define $S_k(z) = 1-G_k(z)$.
The functions $S_k(z)$ satisfy the recurrence
\edef\eqSkformula{\the\eqcount}
$$S_k(z) = S_{k-1}(z) - {z\over S_{k-1}(z)},\adveq$$
with $S_1(z) = 1-z$. 

We define $z^*_k$ to be the smallest positive real root of $S_k(z)$; hence we must have $S_k(z) > 0$ for $z\in[0,z^*_k]$. Furthermore, the formula \refeq{\eqSkformula} implies that $S_k(z)\le S_{k-1}(z)$ for $z\in [0,z^*_k]$, so $z^*_k$ is a nonincreasing sequence of positive numbers.

We now give a lower bound for $z^*_k$ and
show that $S_k(z)$ is small at a point close to this lower bound.

\edef\lemzkstarLB{\the\thmcount}
\proclaim Lemma \advthm. For every $k\ge 1$,
$$z^*_k\ge k-\sqrt{k^2-1}\adveq$$
and
$$S_k\biggl({1\over 2k}\biggr) \le \biggl({1\over 4k}\biggr)^{1/4}.\adveq$$

\proof Squaring both sides of \refeq{\eqSkformula} yields
$$S_k(z)^2-S_{k-1}(z)^2 = -2z + {z^2 \over S_{k-1}(z)^2},\adveq$$
and by telescoping, we have
$$\eqalign{
S_k(z)^2 &= (1-z)^2 + \sum^k_{j=2} \bigl(S_{j}(z)^2 - S_{j-1}(z)\bigr)^2 \cr
& = (1-z)^2 -2(k-1)z + \sum^{k}_{j=2} {z^2 \over S_{j-1}(z)^2} \cr 
& = 1-2kz+z^2 + \sum^{k-1}_{j=1} {z^2 \over S_{j}(z)^2}. 
}\adveq$$
Thus the lower bound $S_k(z)^2 \ge 1-2kz+z^2$ is immediate, and since $1-2kz+z^2$ is positive for $0\le z< k-\sqrt{k^2-1}$, we must have $z^*_k\ge k-\sqrt{k^2-1}$. This proves the first claim.

From the bound just proved, we find that $1/(2k) \in [0,z^*_k]$, so by our earlier observation,
the sequence $S_k\bigl(1/(2k)\bigr)$ is nonincreasing in $k$.
Moreover, we have $S_j(z)\le S_j(0)=1$ for all $j\ge 1$; from this we obtain
$$\eqalign{
S_k\Bigl({1 \over 2k}\Bigr)^2
&= 1 -2k\biggl({1 \over 2k}\biggr) +\biggl({1 \over 2k}\biggr)^2 + \sum^{k-1}_{j=1} {1 \over 4k^2 S_j\bigl(1/(2k)\bigr)^2} \cr 
&\le {1 \over 4k^2} + {1 \over 4k^2}\cdot {k-1\over S_k\bigl(1/(2k)\bigr)^2} \cr 
&\le {1\over 4k S_k\bigl(1/(2k)\bigr)^2}.
}\adveq$$
Solving the inequality gives
$S_k\bigl(1/(2k)\bigr) \le 1/(4k)^{1/4}$, as desired. \slug

Next, we derive an upper bound for $z_k^*$.

\proclaim Lemma \advthm. For all $k\ge 1$,
$$z^*_k\le {1\over 2k\bigl(1-1/(4k)^{1/4}\bigr)}.\adveq$$

\proof
First we observe that $S_k$ is a concave function on $[0,z_k^*]$. This follows from the fact that $G_k$ is a power series with positive coefficients and no constant term, so $S_k = 1-G_k$ has constant term $1$ and all other coefficients negative. Moreover, as $G_k$ is a rational function, it is analytic within its radius of convergence, hence all of its derivatives are well-defined and negative on $[0,z_k^*]$.

Concavity of $S_k$ gives us
$$S_k\biggl({1\over 2k}\biggr)  \ge {1 \over 2kz^*_k} S_k(z^*_k) + \biggl(1 - {1 \over 2kz^*_k}\biggr) S_k(0) = 1 - {1 \over 2kz^*_k},\adveq$$
and applying Lemma~{\lemzkstarLB} now yields
$$\biggl({1\over 4k}\biggr)^{1/4} \ge 1 - {1 \over 2kz^*_k},\adveq$$
and the result follows upon rearranging.\slug

We can now describe the first-order asymptotic behaviour of $G_{n,k}$.

\edef\thmGnk{\the\thmcount}
\proclaim Theorem \advthm. Let $G_{n,k}$ denote the number of planted plane trees with decreasing labels. Then for all $k\ge 1$, there exist real numbers $c_k$ and $\alpha_k$ such that
$$G_{n,k} = c_k{\alpha_k}^n\bigl(1+o(1)\bigr),\adveq$$
as $n\to\infty$, where
$$ 2(k-1)\biggl(1-{1\over \sqrt 2(k-1)^{1/4}}\biggr)
\le \alpha_k \le {1 \over k-1 - \sqrt{(k-1)^2 - 1}}.\adveq$$
In particular, $\alpha_k = 2k+o(k)$ as $k\to\infty$.
\goodbreak

\proof As we did previously, for all $k\ge 1$ let $z_k^*$ denote the solution to
$G_k(z_k^*) = 1$. We know that for $k\ge 2$, $G_k(z)$ has only one pole on the domain $|z|\le z_{k-1}^*$, namely, $z_{k-1}^*$. We claim that this is a simple pole. Since $G_1(z) = z$ and
$$G_k(z) = {G_{k -1}(z) - G_{k-1}(z)^2 + z\over 1-G_{k-1}(z)}\adveq$$
for $k\ge 2$,
we see that $G_k(z)$ is a rational function. Thus
$$\lim_{z\to z_{k-1}^*} (z-z_{k-1}^*)G_k(z)
= {G_{k-1}(z_{k-1}^*)- G_{k-1}(z_{k-1}^*)^2 + z_{k-1}^*\over
-G_{k-1}'(z_{k-1}^*)}
= {z_{k-1}^* \over - G_{k-1}'(z_{k-1}^*)}.\adveq$$
The function $G_{k-1}(z)$ is analytic on the domain $|z|\le z_{k-1}^*$ and its power series has nonnegative coefficients. Since the constant term in the power series $G_{k-1}'(z)$ is
$$[z]G_{k-1}(z)=G_{1,k-1} = k-1 > 0,\adveq$$
the above limit is finite. This
means that $z_{k-1}^*$ is a pole of order $1$ of the function $G_k(z)$.
Letting $\alpha_k = 1/z_{k-1}^*$, by a standard result on the coefficient asymptotics of rational functions, we have
\edef\eqGnkasymptotic{\the\eqcount}
$$G_{n,k} = {c_k {\alpha_k}^n} \bigl(1 + o(1)\bigr),\adveq$$
where $c_k$ is a constant depending on $k$ only.
The claimed bounds on $\alpha_k$ follow from the two previous lemmas.\slug 

We also note the following explicit upper bound.
As $G_k(z_k^*) = 1$, we immediately have
$$
G_{n,k} 
\le \left({ 1 \over z^*_{k}}\right)^n 
= \left(\alpha_{k+1}\right)^n.
$$
This is weaker than~\refeq{\eqGnkasymptotic}, but only slightly so,
as $\alpha_{k+1}/\alpha_k \to 1$ as $k \to \infty$.

\medskip\boldlabel The regular leaning tree of order {\mathbold k}.
Consider the sequence $T_0, T_1, T_2,\ldots$ of trees defined recursively as follows.
Let $T_0$ be a single root node with no children, and for $k\ge 1$, we define $T_k$ to be the tree with a root node having $k$ children whose subtrees are $T_{k-1}, T_{k-2},\ldots, T_0$. We shall call $T_k$ the {\it regular leaning tree of order $k$}. For all $k\ge 0$, the tree $T_k$ has $2^k$ nodes.

The following algorithm gives a mapping from the set of
closed walks of length $2n$ starting at the root of $T_k$ and
$n$-node planted plane trees with root label $k+1$ and decreasing
labels.
\goodbreak

\algbegin Algorithm P (Build planted plane tree). Given a closed walk
$$\sigma = (u_0, u_1, \ldots, u_{2n})$$
of length $2n$ in $T_k$, this
algorithm outputs a planted plane tree labelled with positive integers.
\algstep P1. [Initialize.] Set $i\gets 0$ and initialize the tree $P$ with a root node $v$ and set $\mfield{LABEL}{v} \gets k+1$.
\algstep P2. [Child.] If $u_{i+1}$ is the parent of $u_i$, go to step P3. Otherwise, suppose that $u_{i+1}$ is the root of a subtree isomorphic to $T_j$ for some $j<\mfield{LABEL}{v}$. Append a child with label $j+1$ to $v$ and update $v$ to point to this child. Go to step P4.
\algstep P3. [Parent.] Set $v\gets \mfield{PARENT}{v}$.
\algstep P4. [Loop.] Increment $i$ by one. If $i< 2n$, go to step P2; otherwise terminate the algorithm with the tree $P$ as output.\slug

The root of the output planted plane tree has label $k+1$, and every time we move downward in $T_k$ along the path $\sigma$ and encounter the root of $T_j$ for some $j\le k$, we append a node with label $j+1$. This implies that the output of Algorithm~P is a planted plane tree with decreasing labels and root label $k+1$. We also know that this output tree has $n+1$ nodes, since we start with one root node, and in the walk $\sigma$, half of the $2n$ edges must go down in the tree and half must go up (to return to the root), and we add a node to $P$ when (and only when) going down.

The tree $P$ is built up in depth-first preorder, so it is easy to write an algorithm that recovers the walk from a tree $P$ output by Algorithm~P.

\algbegin Algorithm W (Build walk in $T_k$). Given a planted plane tree $P$ with decreasing labels and root label $k+1$, we build a walk $\sigma$ starting at the root in $T_k$. Suppose we have a function $\mfield{NEXT}{v}$ of getting the node that follows a given node $v\in P$ in a depth-first traversal of $P$ (in this traversal, nodes may be visited multiple times).
\algstep W1. [Initialize.] Let $u_0$ be the root of $T_k$. Set $i\gets 0$ and let $v$ be the root of $P$.
\algstep W2. [Child.] Let $w \gets \mfield{NEXT}{v}$. If $w$ is the parent of $v$, go to step W3. Otherwise, let $u_{i+1}$ be the child
of $u_i$ that is the root of the subtree isomorphic to $T_j$, where
$j=\mfield{LABEL}{w} -1$. Go to step W4.
\algstep W3. [Parent.] Set $u_{i+1} \gets \mfield{PARENT}{u}$.
\algstep W4. [Loop.] Set $v\gets w$ and increment $i$ by one.
(An invariant we have maintained is that $u_i$ is the root of a subtree isomorphic to $T_j$,
where $j = \mfield{LABEL}{v} - 1$. This makes step W2 possible in the next iteration.) If $i< 2n$, go to step W2; otherwise, terminate the algorithm with output $\sigma = (u_0, u_1,\ldots, u_{2n})$.\slug
\goodbreak

An example of a tree $P$ and its corresponding walk in
$T_k$ is illustrated in Fig.~1.
The parallel structures of Algorithms~P and~W make it clear that if Algorithm~P terminates with output $P$ upon being given input $\sigma$, then Algorithm~W returns the walk $\sigma$ upon the input $P$.
We have thus furnished a bijective proof of the following theorem.
\topinsert
$$\epsfbox{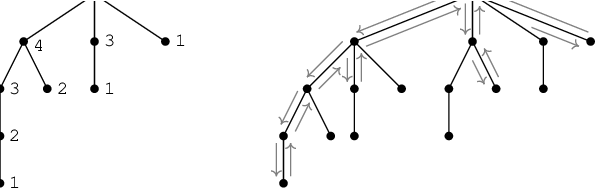}$$
\caption{A planted plane tree with root label $5$ and its corresponding walk in $T_4$.}
\endinsert
\goodbreak

\edef\thmbijection{\the\thmcount}
\proclaim Theorem \advthm. Let $n$ and $k$ be positive integers. The number $W_{2n}(T_k)$ of closed paths in $T_k$ of length $2n$ that begin and end at the root is equal to the number $G_{n+1,k+1} - G_{n+1,k}$
of $(n+1)$-node planted plane trees with decreasing labels and root label equal to $k+1$.\slug

\medskip\boldlabel The top eigenvalue of a regular leaning tree.
Let $\Gamma=(V,E)$ be a graph, let $A = A(\Gamma)$ be its adjacency matrix, and let $\lambda_1(A)$ denote the largest eigenvalue of $A$. By the trace method, we have
$$
\eqalign{
\lambda_1(A) 
&= \lim_{n\to \infty} \Bigl(\tr\bigl(A^{2n}\bigr)\Bigr)^{1/2n} \cr
&= \lim_{n\to\infty} \Bigl( \max_{v\in V} W_{2n}(v,\Gamma)\Bigr)^{1/2n},\cr
&= \lim_{n\to\infty} \Bigl( \min_{v\in V} W_{2n}(v,\Gamma)\Bigr)^{1/2n},\cr
}\adveq
$$
where $W_{2n}(v,\Gamma)$ denotes the number of closed walks of length $2n$ in $\Gamma$ starting at the vertex $v$.

It is well known that if $\Gamma = T$ is a tree and $A$ is its adjacency matrix, then
the largest eigenvalue $\lambda_1(A)$ of $A$ satisfies
$$\sqrt{\Delta} \le \lambda_1(A)\le 2\sqrt{\Delta-1},\adveq$$
where $\Delta$ is the maximum vertex degree of $T$. (The lower bound is trivial and the upper bound is a result of D.~Stevanovi\'c~\ref{\refstevanovic}.)
Theorems~{\thmGnk} and~{\thmbijection} together tell us that
largest eigenvalue of the adjacency matrix of a leaning tree of order $k$ (which has maximum vertex degree $k+1$) does not tend towards either of these bounds as $k\to\infty$.

\edef\lemTkeigenvalue{\the\thmcount}
\proclaim Lemma \advthm. Let $k\ge 1$ and let $A_k$ denote the adjacency matrix of the regular leaning tree $T_k$. Then the
largest eigenvalue $\lambda_1(A_k)$ is $\sqrt{2k + o(k)}$.\slug
\goodbreak

\medskip\boldlabel The Ulam--Harris number.
We can extend the above result to arbitrary trees as follows. We give the root the label $1$, and for any node with label $r$ and $s$ children, we label the children $r+1, r+2,\ldots, r+s$. We define the {\it Ulam--Harris number $\UH(T)$} of a planted plane tree $T$ to be the maximum label in the tree. (There is a standard notion of the {\it Ulam--Harris labelling} of a tree (see, e.g., Section~6 of~\ref{\refjanson}), in which one assigns a vector of positive integers to each node. Our Ulam--Harris number is the maximum sum of coordinates, taken over all Ulam--Harris labels in the tree.) For an {\it unordered} tree $T$, we can let
$\UH(T)$ be the minimum of $\UH(T')$ over all planted plane trees $T'$ obtained from $T$ by choosing orderings for the children of each node.
\topinsert
\vskip10pt
$$\epsfbox{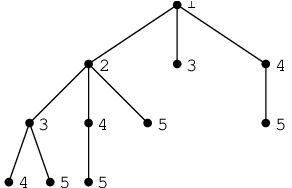}$$
\caption{A planted plane tree with Ulam--Harris number equal to $5$.}
\endinsert

\edef\thmUHTbound{\the\thmcount}
\proclaim Theorem \advthm. Let $T$ be a rooted unordered tree, and let $A$ be the adjacency matrix of $T$. Then $\lambda_1(A)\le \lambda_1\bigl(A_{\UHsmall(T)-1}\bigr)$,
where $A_{\UHsmall(T)-1}$ is the adjacency matrix of the regular leaning tree of order $\UH(T)-1$.

\proof Consider $T$ as a labelled planted plane tree, with maximum label $\UH(T)$. Replace each label $j$ with $\UH(T) - j + 1$. Then the root has label $\UH(T)$, all labels are positive and decreasing both as one descends in the tree, and as one goes from left to right among siblings. Furthermore, the new label of each node is strictly greater than the number of children it has, so we can embed $T$ into the regular leaning tree $T_{\UHsmall(T)-1}$, and the claim follows.\slug

Let $T$ be a tree with maximum vertex degree $\Delta$ and adjacency matrix $A$.
Since a node and its children must all have different labels, we must have $\UH(T) \ge \Delta$; the regular leaning tree of order $k$ attains this bound with equality. Lemma~{\lemTkeigenvalue} and Theorem~{\thmUHTbound} together supply an upper bound of roughly $\sqrt{2\UH(T)-2}$ on $\lambda_1(A)$. In the scenario where $\UH(T)+1 < 2\Delta$ (as in Fig.~2, for instance), this improves upon Stevanovi\'c's bound of $2\sqrt{\Delta-1}$.

\medskip\boldlabel Related work.
In this paper we have considered trees with decreasing labels, but in the literature it has been somewhat more common to count trees with {\it increasing} labels (by reflecting the values of the nodes, there are exactly as many increasing trees with maximum label $k$ as there are decreasing ones). There is a classical bijection between permutations of $\{1,\ldots n\}$ and increasing binary trees on $n$ nodes using each label from $1$ through $n$ exactly once. A 2020 paper of O.~Bodini, A.~Genitrini, B.~Gittenberger, and S.Wagner~\ref{\refbodinione} relaxed the stipulation that all the labels $1$ through $n$ must appear, requiring instead that for some $1\le k\le n$, all the labels $1$ through $k$ must appear (labels must still strictly increase down the tree). A later paper of O.~Bodini, A.~Genitrini, M.~Naima, and A.~Singh describes a more general approach that applies to more families of trees~\ref{\refbodinitwo}.

In 1987, J.~Blieberger counted the number of Motzkin trees with labels that increase, but not necessarily strictly~\ref{\refblieberger}.
R.~Kemp showed in 1993 that the planted plane trees labelled similarly are in bijection with monotonically extended binary trees~\ref{\refkemp}.
In 2011, S.~Janson, M.~Kuba, and A.~Panholzer drew a link between generalized Stirling permutations and $k$-ary trees with (strictly) increasing labels~\ref{\refjansonkuba}. A generalization of labelled trees, in which nodes can receive multiple labels, was introduced in 2016 by M.~Kuba and A.~Panholzer, and this generalization is shown to imply various hook-length formulas for trees~\ref{\refkubatwo}.

Other authors have considered the average shape~\ref{\refkirschenhofer} and the degree distribution~\ref{\refkubaone} of various tree families with increasing labels. The expectation and variance of the size of the ancestor tree as well as the Steiner distance of increasingly labelled trees were determined by K.~Morris in 2004~\ref{\refmorris}.

\medskip\boldlabel Acknowledgements.
We thank Ari Blondal, G\'abor Lugosi, and Steve Melczer for stimulating discussions and helpful suggestions. The second and third authors are funded by the Natural Sciences and Engineering Research Council of Canada. 

\medskip
\section References

\def\beginref{\noindent}
\def\endref{\medskip}
\begingroup\frenchspacing\tolerance=2200

\beginref
\parindent=20pt\item{\ref{\refblieberger}}
Johann Blieberger,
``Monotonically labelled Motzkin trees,''
{\sl Discrete Applied Mathematics}
{\bf 18}
(1987),
9--24.
\endref

\beginref
\parindent=20pt\item{\ref{\refbodinione}}
Olivier Bodini,
Antoine Genitrini,
Bernhard Gittenberger,
and Stephan Wagner,
``On the number of increasing trees with label repetitions,''
{\sl Discrete Mathematics}
{\bf 343}
(2020),
article no.~111722.
\endref

\beginref
\parindent=20pt\item{\ref{\refbodinitwo}}
Olivier Bodini,
Antoine Genitrini,
Mehdi Naima,
and Alexandros Singh,
``Families of monotonic trees: Combinatorial enumeration and asymptotics,''
in {\sl 15th International Computer Science Symposium in Russia} (2020), 155--168.
\endref

\beginref
\parindent=20pt\item{\ref{\refanalyticcombinatorics}}
Philippe Flajolet
and Robert Sedgewick,
{\sl Analytic Combinatorics}
(Cambridge:
Cambridge University Press,
2009).
\endref

\beginref
\parindent=20pt\item{\ref{\refjanson}}
Svante Janson,
``Simply generated trees,  conditioned Galton--Watson trees,  random allocations and condensation,''
{\sl Probability Surveys}
{\bf 9}
(2012),
103--252.
\endref

\beginref
\parindent=20pt\item{\ref{\refjansonkuba}}
Svante Janson,
Markus Kuba,
and Alois Panholzer,
``Generalized Stirling permutations, families of increasing trees and urn models,''
{\sl Journal of Combinatorial Theory, Series A}
{\bf 118}
(2011),
94--114.
\endref

\beginref
\parindent=20pt\item{\ref{\refkemp}}
Rainer Kemp,
``Monotonically labelled ordered trees and multidimensional binary trees,''
in {\sl International Symposium on Fundamentals of Computation Theory} (1993),
329--341.
\endref

\beginref
\parindent=20pt\item{\ref{\refkirschenhofer}}
Peter Kirschenhofer,
``On the average shape of monotonically labelled tree structures,'' {\sl Discrete Applied Mathematics}
{\bf 7}
(1984),
161--181. 
\endref

\beginref
\parindent=20pt\item{\ref{\refkubaone}}
Markus Kuba
and Alois Panholzer,
``On the degree distribution of the nodes in increasing trees,''
{\sl Journal of Combinatorial Theory, Series A}
{\bf 114}
(2007),
597--618.
\endref

\beginref
\parindent=20pt\item{\ref{\refkubatwo}}
Markus Kuba
and Alois Panholzer,
``Combinatorial families of multilabelled increasing trees and hook-length formulas,''
{\sl Discrete Mathematics}
{\bf 339}
(2016)
227--254.
\endref

\beginref
\parindent=20pt\item{\ref{\refmorris}}
Katherine Morris,
``On parameters in monotonically labelled trees,''
in {\sl Mathematics and Computer Science {\mc III}}
(2004),
261--263.
\endref

\beginref
\parindent=20pt\item{\ref{\refpemantle}}
Robin Pemantle,
Mark Curtis Wilson,
and Stephen Melczer,
{\sl Analytic Combinatorics in Several Variables}
(Cambridge:
Cambridge Studies in Advanced Mathematics,
2024).
\endref

\beginref
\parindent=20pt\item{\ref{\refstevanovic}}
Dragan Stevanovi\'c,
``Bounding the largest eigenvalue of trees in terms of the largest vertex degree,''
{\sl Linear Algebra and its Applications}
{\bf 60}
(2003),
35--42.
\endref

\endgroup

\bye